\documentclass{amsart}

\usepackage{epsfig}
\usepackage{graphicx}

\newtheorem{theorem}{Theorem}[section]

\theoremstyle{definition}

\newtheorem{example}[theorem]{Example}

\theoremstyle{remark}
\newtheorem{remark}[theorem]{Remark}

\numberwithin{equation}{section}

\begin{document}

\title{Chebyshev Series Expansion of Inverse Polynomials}

\author{Richard J. Mathar}
\address{Sterrewacht, Universiteit Leiden, Postbus 9513, 2300 RA Leiden, The Netherlands}
\email{mathar@strw.leidenuniv.nl}

\subjclass[2000]{Primary 33C45, 42C20; Secondary 41A50, 41A10; Tertiary 65D15}

\date{\today}

\keywords{Chebyshev series, orthogonal polynomials, approximation}

\begin{abstract}
An inverse polynomial has a Chebyshev series expansion
$$1/\sum_{j=0}^kb_jT_j(x)=\mathop{{\sum}'}_{n=0}^\infty a_nT_n(x)$$
if the polynomial has no roots in $[-1,1]$. If the inverse polynomial is decomposed into
partial fractions, the $a_n$ are linear combinations of simple functions
of the polynomial roots. Also, if the first $k$ of the coefficients $a_n$
are known, the others become linear combinations of these with expansion
coefficients derived recursively from the $b_j$'s.
On a closely related theme,
finding a polynomial with minimum {\em relative\/} error towards
a given $f(x)$
is approximately equivalent to
finding the $b_j$ in $f(x)/\sum_0^kb_jT_j(x)=1+\sum_{k+1}^\infty a_nT_n(x)$,
and may be handled with a Newton method providing the Chebyshev
expansion of $f(x)$ is known.
\end{abstract}

\maketitle

\section{Introduction and Scope}
The Chebyshev polynomials $T_n(x)$ are even or odd functions of $x$
defined as \cite[(22.3.6)]{AS}\cite[(3.6)]{Bulychev}
\begin{equation}
T_0(x)=1,\quad T_n(x)=\frac{n}{2}\sum_{m=0}^{\lfloor n/2\rfloor}(-)^m\frac{(n-m-1)!}{m!(n-2m)!}(2x)^{n-2m},
\quad n=1,2,3\ldots
\label{eq.Tnofx}
\end{equation}
where the Gauss bracket $\lfloor .\rfloor$ denotes the largest integer not greater
than the number it embraces. The reverse formula is \cite[p.\ 412]{Cody}\cite[p.\ 52]{Fox}\cite{FraserJACM12}
\begin{equation}
x^n=2^{1-n}{\mathop{{\sum}'}_{\genfrac{}{}{0pt}{}{j=0}{n-j \text{even}}}^n}\binom{n}{(n-j)/2} T_j(x)
\label{eq.xnofT}
\end{equation}
where the prime at the sum symbol means the first term (at $j=0$ and even $n$) is to be halved.
The polynomials are orthogonal over the interval $[-1,1]$
with weight function $1/\sqrt{1-x^2}$ \cite[(22.2.4)]{AS}\cite[(4.2)]{Cody}\cite{BarrioAMC150}
\begin{equation}
\int_{-1}^1 T_n(x)T_m(x)\frac{dx}{\sqrt{1-x^2}}=\left\{
\begin{array}{c@{,\quad}c}
\pi & n=m=0,\\
\pi/2 & n=m\neq 0,\\
0 & n\neq m.\\
\end{array}
\right.
\end{equation}
The product rule is \cite[(22.7.24)]{AS}\cite[(A.1)]{Hasegawa1991}\cite[(2.10)]{FraserJACM12}
\begin{equation}
T_n(x)T_m(x)=\frac{1}{2}\left(T_{|m-n|}(x)+T_{m+n}(x)\right).
\label{eq.Tnmsum}
\end{equation}
The indefinite integral is \cite[(4.8)]{Cody}\cite[p.\ 54]{Fox}\cite[(2.12)]{FraserJACM12}
\begin{equation}
\int T_n(x)dx=\left\{ \begin{array}{ll}
T_1(x), & n=0,\\
\frac{1}{4}T_2(x), & n=1,\\
\frac{1}{2}\left( \frac{T_{n+1}(x)}{n+1}-\frac{T_{n-1}(x)}{n-1} \right), & n>1,\\
\end{array}\right.
\label{eq.Tintx}
\end{equation}
which correlates to the derivative
\begin{equation}
\frac{d}{dx}T_n(x) = 2n
\mathop{{\sum}'}_{\genfrac{}{}{0pt}{}{l=0}{n-l\, \text{odd}}}
^{n-1}T_l(x).
\label{eq.Tderi}
\end{equation}

The expansion of an inverse polynomial of degree $k$ in a power series
is \cite[(3.6.16)]{AS}
\begin{equation}
\frac{1}{\sum_{j=0}^k d_jx^j}=\sum_{n=0}^\infty c_nx^n,
\label{eq.powser}
\end{equation}
with recursively accessible \cite[0.313]{GR}
\begin{equation}
c_n=-\frac{1}{d_0}\sum_{\genfrac{}{}{0pt}{}{j=0}{n-j\le k}}^{n-1}d_{n-j}c_j,
\quad c_0=\frac{1}{d_0},\quad n\ge 1.
\end{equation}

The topic of this script is the equivalent arithmetic expansion of the inverse
polynomial in a Chebyshev series,
\begin{equation}
\frac{1}{\sum_{j=0}^k d_jx^j}=
\frac{1}{\sum_{j=0}^k b_jT_j(x)}=
{\mathop{{\sum}'}_{n=0}^\infty} a_nT_n(x),
\label{eq.anDef}
\end{equation}
i.e., computation of the coefficients
\begin{equation}
a_n=\frac{2}{\pi}\int_{-1}^1 \frac{T_n(x)}{\sum_{j=0}^k b_jT_j(x)}\frac{dx}{\sqrt{1-x^2}}
\label{eq.anofT}
\end{equation}
given the sets \{$b_j$\} or \{$d_j$\} that define the original function.
Both sets are related via \cite[(3)]{Schonfelder}\cite[(37)]{ScratonMathComp50}
and with \eqref{eq.Tnofx} via
\begin{equation}
d_l=\frac{2^{l-1}}{l!}\sum_{\genfrac{}{}{0pt}{}{j=0}{j-l \text{even}}}^k (-)^{(j-l)/2}j
\frac{(\frac{j+l}{2}-1)!}{(\frac{j-l}{2})!}b_j,\quad l=0,\ldots,k.
\label{eq.dlofbl}
\end{equation}

The expansion \eqref{eq.anDef} exists if the inverse polynomial is bound in 
the interval $[-1,1]$,
i.e., if $\sum d_jx^j$ has no roots in $[-1,1]$.

Characteristic generic methods of evaluating \eqref{eq.anofT} are not reviewed here:
(i) Fourier transform methods \cite[(4.7)]{Cody}\cite{Clenshaw,Fosdick,Fettis},
(ii)
sampling with Gauss-type quadratures \cite[(25.4.38)]{AS}\cite{Perez,SmithHV,Kzaz} which effectively
means using an implicit intermediate interpolatory polynomial
to represent $1/\sum_{j=0}^kb_jT_j(x)$,
(iii)
approximation by truncation of \eqref{eq.powser}, then insertion of \eqref{eq.xnofT},
(iv)
using the near-minimax properties of the Chebyshev series \cite{Murnaghan,MasonNumAlg38}.
\begin{remark}
The Fourier-Chebyshev series \cite{Holub,RivlinNumMath4}\cite[(3.4.1f)]{Rivlin}
\begin{equation}
\frac{T_n(z)-tT_{|n-m|}(z)}{1-2tT_m(z)+t^2}=\sum_{k=0}^\infty T_{km+n}(z)t^k
\end{equation}
provides special cases of polynomials with particularly simple expansions.
\end{remark}
\begin{remark}
Insertion of $n=1$ in (\ref{eq.Tnmsum}) shows that the coefficients of
\begin{equation}
f(x)= \mathop{{\sum}'}_{n=0}^\infty f_nT_n(x)
\end{equation}
and
\begin{equation}
\frac{f(x)}{x}= \mathop{{\sum}'}_{n=0}^\infty g_nT_n(x)
\end{equation}
are related as
\begin{equation}
f_0=g_1,\quad 2f_{n-1}=g_{n-2}+g_n,\quad n\ge 2.
\label{eq.fxx}
\end{equation}
\end{remark}

Chapter \ref{pfrac.sec} explains how the $a_n$ of \eqref{eq.anofT}
could be computed supposed the inverse polynomial has been decomposed into
partial fractions. Chapter \ref{recu.sec} provides a recursive algorithm
to derive high-indexed $a_n$ ($n\ge k$) supposed the low-indexed $a_n$
are given by other means. Chapter \ref{trunc.sec} touches on a
(standard) integral-free method to compute approximate  low-indexed $a_n$,
and Chapter \ref{rel.sec} deals with a specific inverse problem
---which is finding the $b_j$ from partially known $a_n$---related to
polynomial approximants with minimum relative error.

\section{The Case of Known Partial Fractions}\label{pfrac.sec}
The straight way of computing the Chebyshev series uses the
decomposition of $1/\sum d_jx^j$ into partial fractions \cite[2.102]{GR},
which reduces \eqref{eq.anDef} to the calculation of the $a_{n,s}$ in
\begin{equation}
\frac{1}{(z-x)^s}\equiv\mathop{{\sum}'}_{n=0}^\infty a_{n,s}(z)T_n(x),
\label{eq.ansDef}
\end{equation}
where $z$ is a root of the polynomial,
\begin{equation}
\sum_{j=0}^kd_j z^j=0.
\end{equation}

Sign flips of $z$ and $x$ in \eqref{eq.ansDef} show that
\begin{equation}
a_{n,s}(-z)=(-)^{n+s}a_{n,s}(z).
\end{equation}
The case of $s=1$ has been evaluated earlier \cite[(A.6)]{Hasegawa1983}\cite{RivlinNumMath4,Sen}
based on \cite[(22.9.9)]{AS}\cite[(18)]{PolP1},
\begin{equation}
a_{n,1}(z)=\frac{2}{(z^2-1)^{1/2}}\frac{1}{w^n},
\quad w\equiv z+(z^2-1)^{1/2},
\quad z\notin [-1,1].
\label{eq.Haseg}
\end{equation}
The branch cuts of $(z^2-1)^{1/2}$ must be chosen such that $|w|>1$.

\begin{example}
\begin{equation}
\frac{1}{1+x^2}=\frac{i}{2}\frac{1}{i-x}-\frac{i}{2}\frac{1}{-i-x}
\end{equation}
consists of two terms,
\begin{equation}
a_{n,1}(i)=-\frac{\sqrt{2}i^{1-n}}{(1+\sqrt{2})^n},\quad
a_{n,1}(-i)=(-)^n\frac{\sqrt{2}i^{1-n}}{(1+\sqrt{2})^n},
\end{equation}
which recombine with the two factors $i/2$ and $-i/2$ to \cite[(3.4.1a)]{Rivlin}
\begin{equation}
\frac{1}{1+x^2}=\sqrt{2}\mathop{{\sum}'}_{n=0,2,4,6,\ldots}\frac{(-)^{n/2}}{(1+\sqrt{2})^n}T_n(x).
\label{eq.1xx}
\end{equation}
\end{example}

\begin{remark}
The shifted Chebyshev polynomials
$T^*(x)\equiv T(2x-1)$ are orthogonal over $[0,1]$ with weight $1/\sqrt{x(1-x)}$
\cite[(22.2.8)]{AS}\cite{RababahCMAM3}. From (\ref{eq.ansDef}) we get
\begin{equation}
\frac{1}{(z-x)^s}=2^s\mathop{{\sum}'}_{n=0}^\infty a_{n,s}(2z-1)T_n^*(x),
\label{eq.ansstarDef}
\end{equation}
and (\ref{eq.Tintx}) becomes
\begin{equation}
\int T_n^*(x)dx=\left\{ \begin{array}{ll}
\frac{1}{2}T_1^*(x), & n=0,\\
\frac{1}{8}T_2^*(x), & n=1,\\
\frac{1}{4}\left( \frac{T_{n+1}^*(x)}{n+1}-\frac{T_{n-1}^*(x)}{n-1} \right), & n>1.\\
\end{array}\right.
\label{eq.Tstarintx}
\end{equation}
\end{remark}

\begin{example}
An example of $s=1$, $z=-1$ in (\ref{eq.ansstarDef}) is
\begin{equation}
\frac{1}{1+x}=-\frac{1}{-1-x}=-2\mathop{{\sum}'}_{n=0}^\infty a_{n,1}(-3)T_n^*(x),
\label{eq.T1x}
\end{equation}
where
\begin{equation}
a_{n,1}(-3)=\frac{(-)^{n+1}}{\sqrt{2}(3+2\sqrt{2})^n}
\label{eq.an-3}
\end{equation}
according to (\ref{eq.Haseg}).
\end{example}

Higher second indices $s$ of the $a_{n,s}$ are obtained from \eqref{eq.ansDef} by repeated
derivation w.r.t.\ $z$,
\begin{equation}
(-)^ss!\frac{1}{(z-x)^{s+1}}={\mathop{{\sum}'}_{n=0}^\infty} \left(\frac{\partial}{\partial z}
\right)^s a_{n,1}(z)T_n(x),
\end{equation}
via \cite[0.432.1]{GR},
\begin{eqnarray}
a_{0,s+1}(z)&=&\frac{2}{s!}(-)^s\left(\frac{\partial}{\partial z}\right)^s
\frac{1}{(z^2-1)^{1/2}} \nonumber \\
&=& 2\sum_{l=0}^{\lfloor s/2\rfloor}\frac{(-)^l}{l!(s-2l)!}\left(\frac{1}{2}\right)_{s-l}
\frac{(2z)^{s-2l}}{(z^2-1)^{\frac{1}{2}+s-l}},\quad s\ge 0, \label{eq.a2ndRec}
\end{eqnarray}
with Pochhammer's Symbol defined as \cite[(6.1.22)]{AS}
\begin{equation}
(\alpha)_k\equiv \alpha(\alpha+1)(\alpha+2)\cdots (\alpha+k-1)=\Gamma(\alpha+k)/\Gamma(\alpha),
\quad (\alpha)_0=1.
\end{equation}
The formula
\begin{eqnarray}
(s-1)\int \frac{dx}{(z-x)^s}&=&\frac{1}{(z-x)^{s-1}}+A_s \nonumber \\
&=&
(s-1) \mathop{{\sum}'}_{n=0}^\infty a_{n,s}\int T_n(x)dx
= \mathop{{\sum}'}_{n=0}^\infty a_{n,s-1}T_n(x)+A_s,
\quad s\ge 2 \nonumber
\end{eqnarray}
in conjunction with the method quoted by Cody \cite[(4.8)]{Cody}\cite[(25)]{Mason1999} yields
\begin{equation}
a_{n+1,s}(z)=a_{n-1,s}(z)-\frac{2n}{s-1}a_{n,s-1}(z),
\quad n\ge 1,\quad s\ge 2. \label{eq.a1stRec}
\end{equation}
One needs \eqref{eq.a2ndRec} and
\begin{eqnarray}
a_{1,s+1}(z)&=&\frac{2}{\pi}\int_{-1}^1 \frac{T_1(x)}{(z-x)^{s+1}}\frac{dx}{\sqrt{1-x^2}} \nonumber \\
&=&
-a_{0,s}(z)+za_{0,s+1}(z) \label{eq.powRec}
\end{eqnarray}
to start  the recurrence \eqref{eq.a1stRec} and to obtain all coefficients in
\eqref{eq.ansDef} for a particular $z$. Closed form expressions for
solving these recurrences in terms of Legendre Polynomials of
$z/\sqrt{z^2-1}=(w^2+1)/(w^2-1)$ have been
given by Elliott \cite{ElliottMathComp18}.

\begin{remark}
\eqref{eq.powRec} may be generalized to
\begin{equation}
\frac{2}{\pi}\int_{-1}^1 \frac{x^l}{(z-x)^n}\frac{dx}{\sqrt{1-x^2}} 
=\sum_{m=0}^l (-)^m\binom{l}{m} z^{l-m}a_{0,n-m},\quad l<n.
\end{equation}
and with \eqref{eq.xnofT} and \eqref{eq.Tnmsum} to
\begin{equation}
\frac{2}{\pi}\int_{-1}^1 \frac{x^l}{(z-x)^n}\frac{T_s(x)}{\sqrt{1-x^2}} dx
= \frac{1}{2^l} \mathop{{\sum}'}_{\genfrac{}{}{0pt}{}{i=0}{l-i \text{even}}}^l \binom{l}{\frac{l-i}{2}}
\left[ a_{|i-s|,n}+a_{i+s,n}\right].
\end{equation}
\end{remark}

\begin{example}
An example of degree $k=3$
is
\begin{eqnarray}
\frac{1}{(4-x)^2(5+x)}
&=&\frac{1}{78\frac{1}{2}T_0(x)-23\frac{1}{4}T_1(x)-1\frac{1}{2}T_2(x)+\frac{1}{4}T_3(x)} \label{eq.exa2}\\
&=&\frac{1}{9}\frac{1}{(4-x)^2}+\frac{1}{81}\frac{1}{(4-x)}-\frac{1}{81}\frac{1}{(-5-x)}.
\label{eq.exa1}
\end{eqnarray}
The root at $z=4$ yields
\begin{equation}
a_{0,1}(4)=2/\sqrt{15}\approx 0.5164
\end{equation}
from \eqref{eq.Haseg} and
\begin{equation}
a_{0,2}(4)=2\cdot\frac{1}{2}\cdot\frac{2\cdot 4}{\sqrt{15}^3}\approx 0.1377
\end{equation}
from \eqref{eq.a2ndRec}. The root at $z=-5$ yields
\begin{equation}
a_{0,1}(-5)=2/(-\sqrt{24})\approx -0.4082
\end{equation}
from \eqref{eq.Haseg}. The combined total in \eqref{eq.exa1} is
\begin{equation}
a_0=\frac{2}{\pi}\int_{-1}^1\frac{1}{(4-x)^2(5+x)}\frac{dx}{\sqrt{1-x^2}}
\approx \frac{1}{9}\cdot 0.1377+\frac{1}{81}\cdot0.5164-\frac{1}{81}\cdot (-0.4082)
\approx 0.0267.
\end{equation}
\end{example}
\begin{example}
A case of $k=\infty$ is \cite[1.421.2]{GR} 
\begin{equation}
\frac{\tanh(\pi x/2)}{x}=\frac{4}{\pi}\sum_{m=1}^\infty\frac{i}{2(2m-1)}\left[
\frac{1}{i(2m-1)-x}-\frac{1}{-i(2m-1)-x}\right].
\end{equation}
The roots at $z=\pm i(2m-1)$ yield
\begin{equation}
a_{0,1}(z)=2/\left(\pm i\sqrt{4m^2-4m+2}\right),
\end{equation}
and the combined total is
\begin{equation}
a_0=\frac{2}{\pi}\int_{-1}^1\frac{\tanh(\pi x/2)}{x}\frac{dx}{\sqrt{1-x^2}}
= \frac{4}{\pi}\sum_{m=1}^\infty \frac{1}{(2m-1)\sqrt{m^2-m+1/2}}
\approx 2.38.
\end{equation}
\end{example}

\begin{remark}
From \eqref{eq.Haseg}
\begin{equation}
\frac{\partial a_{n,1}(z)}{\partial z}=
-a_{n,1}\left[\frac{z}{z^2-1}+\frac{n}{(z^2-1)^{1/2}}\right],
\end{equation}
so the (linear) propagation of the absolute relative error in the
root $z$ to the error in the coefficient $a_{n,1}$ is
\begin{equation}
\left| \frac{\Delta a_{n,1}}{a_{n,1}}\right|=
\left| \frac{\Delta z}{z}\right|\cdot\left|\frac{z^2}{z^2-1}+\frac{nz}{(z^2-1)^{1/2}} \right|.
\end{equation}
\end{remark}

\begin{remark}
An associated factorization $\sum_{j=0}^k b_jT_j(x)\propto \prod_{m=1}^l (z_m-x)^{s_m}$,
with $l$ different roots of multiplicities $s_m$, decomposes the
square root of the polynomial
into a $l$-fold product of series of the prototypical forms
\begin{eqnarray}
\sqrt{z-x}&=& \mathop{{\sum}'}_{n=0}^\infty q_n(z)T_n(x),\quad s_m=1, \\
z-x&=& zT_0(x)-T_1(x),\quad s_m=2,
\end{eqnarray}
where \cite[2.576.2]{GR}
\begin{equation}
q_0(z)=\frac{2}{\pi}\int_{0}^\pi dt\sqrt{z-\cos t}=\frac{4}{\pi}\sqrt{1+z}E(\frac{2}{1+z})
\end{equation}
is related to Complete Elliptic Integrals of the Second Kind $E$ in the
notation of \cite[(17.3.4)]{AS}.
The $q_n(z)$ with $n\ge 1$ follow recursively using \cite[(17.1.4)]{AS}. In particular,
one may expand $T_n(x)$ in terms of $P_n^{(0,-1/2)}(x)$ with \cite[(1.4)]{Fields}
to obtain
\begin{equation}
q_n(1)=\frac{2^{5/2}}{\pi}\sum_{l=0}^n\frac{(-n)_l(n)_l}{(3/2)_l(1/2)_l},
\quad n=0,1,2,\ldots
\end{equation}
for the Chebyshev coefficients of $\sqrt{1-x}$.
See \cite{SawadaLNCS2517} for an application.
\end{remark}

\section{Recurrence of Expansion Coefficients}\label{recu.sec}
The $T_n$ in \eqref{eq.anofT} may be decomposed into a unique product of
a polynomial by the denominator plus a remainder of polynomial degree less than $k$.
[The dependence on $x$ is omitted at all $T_n(x)$ for brevity.]
\begin{eqnarray}
T_n&=&(d_0^{(n)}T_0+d_1^{(n)}T_1+\cdots+d_{n-k}^{(n)}T_{n-k})
(b_0T_0+b_1T_1+\cdots+b_kT_k) \nonumber\\
&& +\frac{c_0^{(n)}}{2}T_0+c_1^{(n)}T_1+c_2^{(n)}T_2+\cdots+c_{k-1}^{(n)}T_{k-1}.
\label{eq.Tndiv}
\end{eqnarray}
Expansion with \eqref{eq.Tnmsum} yields a system of linear equations for the vector of the
unknowns $d_j^{(n)}$ and $c_j^{(n)}$:
\begin{equation}
\left(
\begin{array}{ccccc@{|}ccccc}
1 & 0 & \ldots & \ldots & 0 & 2b_0 & b_1  & b_2 & b_3 & \ldots\\
0 & 1 & 0 & \ldots & 0 & b_1 & b_0+\frac{b_2}{2} & \frac{b_1+b_3}{2} & \frac{b_2+b_4}{2} & \ldots \\
\vdots & 0 & \ddots & \ddots & \vdots & b_2 & \frac{b_1+b_3}{2}  & b_0+\frac{b_4}{2} &\frac{b_1+b_5}{2}& \ldots  \\
\vdots & \vdots & \ddots & 1 & 0 & \vdots  & \vdots  & \vdots & \vdots & \vdots \\
0 & \ldots & \ldots & 0 & 1 & b_{k-1} & \frac{b_{k-2}+b_k}{2}  & \frac{b_{k-3}}{2} & \ldots & \ldots  \\
\hline
0 & \ldots & \ldots & \ldots & 0 & b_k & \frac{b_{k-1}}{2} & \frac{b_{k-2}}{2} & \frac{b_{k-3}}{2} & \ldots \\
\vdots & \ldots & \ldots & \ldots & \vdots & 0 & \frac{b_k}{2} & \frac{b_{k-1}}{2} & \ldots & \ldots \\
\vdots & \ldots & \ldots & \ldots & \vdots & \vdots  & 0 & \\
\vdots & \vdots & \vdots & \vdots & \vdots & \vdots & \vdots & \ddots & & \frac{b_{k-1}}{2} \\
0 & \ldots & \ldots & \ldots & 0 & 0 & \ldots & \ldots & 0 & \frac{b_k}{2} \\
\end{array}
\right)
\cdot
\left(
\begin{array}{c}
c_0^{(n)} \\
c_1 ^{(n)}\\
\vdots \\
\vdots \\
c_{k-1} ^{(n)}\\
\hline
d_0 ^{(n)}\\
d_1 ^{(n)}\\
\vdots \\
d_{n-k-1}^{(n)} \\
d_{n-k}^{(n)}
\end{array}
\right)
=
\left(
\begin{array}{c}
0 \\
\vdots \\
\vdots \\
\vdots \\
0 \\
\hline
0 \\
\vdots \\
\vdots \\
0 \\
1
\end{array}
\right)
\label{eq.cdRec}
\end{equation}
The $(n+1)\times (n+1)$ coefficient matrix $A_{r,c}$
(row index $r$ and column index $c$ from 0 to $n$) is an
upper triangular matrix.
It hosts a $k\times k$
unit matrix in the upper left corner, and is symmetric w.r.t.\ the minor
diagonal that stretches from $A_{0,k}$ to $A_{n-k,n}$:
\begin{equation}
A_{r,c}=\delta_{r,c},\quad 0\le c\le k-1.
\end{equation}
\begin{equation}
A_{r,k+c}=A_{c,k+r}=\left\{
\begin{array}{c@{,\quad}c}
2b_0 & r=c=0 \\
b_c & r=0,\quad 1\le c\le k \\
\frac{1}{2}(b_{|r-c|}+b_{r+c}) & r\neq c 
\\
b_0+b_{2c}/2 & r=c=1,2,\ldots,n-k \\
\end{array}
\right.
\end{equation}
This works with the auxiliary definition
\begin{equation}
b_i=0, \quad i>k\quad \text{or}\quad i<0.
\label{eq.bbeyond}
\end{equation}
Insertion of \eqref{eq.Tndiv} into \eqref{eq.anofT} yields
\begin{equation}
a_n=2d_0^{(n)}+ \mathop{{\sum}'}_{i=0}^{k-1} c_i^{(n)}a_i, \quad n\ge k,
\label{eq.anredu}
\end{equation}
which means that entire sequence $a_n$ can be generated recursively from 
its first $k$ terms, if the $d_0^{(n)}$ and $c_i^{(n)}$
are generated at the same time via \eqref{eq.cdRec} or an equivalent method.
Iterated full solution of \eqref{eq.cdRec} can be avoided
through recursive generation of the set \{$d_i^{(n+1)},c_i^{(n+1)}$\}
from
\{$d_i^{(n)},c_i^{(n)}$\}
and
\{$d_i^{(n-1)},c_i^{(n-1)}$\} as follows:
\begin{eqnarray}
d_0^{(n+1)}&=& d_1^{(n)}+\frac{c_{k-1}^{(n)}}{b_k}-d_0^{(n-1)} , \label{eq.dcRec1st}\\
d_1^{(n+1)}&=& 2d_0^{(n)}+d_2^{(n)}-d_1^{(n-1)} ,\\
d_j^{(n+1)}&=& d_{j-1}^{(n)}+d_{j+1}^{(n)}-d_j^{(n-1)},\quad j=2,3,\ldots ,n-k+1.\\
\frac{c_0^{(n+1)}}{2}&=& c_1^{(n)}-\frac{b_0c_{k-1}^{(n)}}{b_k}-\frac{c_0^{(n-1)}}{2} ,\\
c_j^{(n+1)}&=& c_{j-1}^{(n)}+c_{j+1}^{(n)}-\frac{b_jc_{k-1}^{(n)}}{b_k}-c_j^{(n-1)} , \quad j=1,2,\ldots,k-1, \label{eq.dcReclast}
\end{eqnarray}
where the auxiliary definitions 
\begin{equation}
c_j^{(n)}=0,\quad j\ge k,\quad \text{or}\quad j<0,
\end{equation}
\begin{equation}
d_j^{(n)}=0, \quad j>n-k,\quad \text{or}\quad j<0,
\end{equation}
are made to condense the notation.
\begin{proof}
Multiply \eqref{eq.Tndiv} by $2T_1$ and use \eqref{eq.Tnmsum} as
\begin{eqnarray}
2T_1\sum_{j=0}^{n-k}d_j^{(n)}T_j&=&d_1^{(n)}T_0+(2d_0^{(n)}+d_2^{(n)})T_1 \nonumber \\
&& +
\sum_{j=2}^{n-k-1}(d_{j-1}^{(n)}+d_{j+1}^{(n)})T_j+d_{n-k-1}^{(n)}T_{n-k}+d_{n-k}^{(n)}T_{n-k+1},
\end{eqnarray}
\begin{equation}
2T_1\mathop{{\sum}'}_{j=0}^{k-1}c_j^{(n)}T_j=c_1^{(n)}T_0+
\sum_{j=1}^{k-2}(c_{j-1}^{(n)}+c_{j+1}^{(n)})T_j+c_{k-2}^{(n)}T_{k-1}+c_{k-1}^{(n)}T_k.
\end{equation}
Rewrite the last term in the previous equation
\begin{equation}
c_{k-1}^{(n)}T_k=\frac{c_{k-1}^{(n)}}{b_k}\sum_{j=0}^kb_jT_j
-\frac{c_{k-1}^{(n)}}{b_k}b_0T_0-\ldots
-\frac{c_{k-1}^{(n)}}{b_k}b_{k-1}T_{k-1}.
\end{equation}
Construct
\begin{eqnarray}
2T_1T_n&=&\big[
(d_1^{(n)}+\frac{c_{k-1}^{(n)}}{b_k})T_0+(2d_0^{(n)}+d_2^{(n)})T_1 \nonumber \\
&& \quad +
\sum_{j=2}^{n-k-1}(d_{j-1}^{(n)}+d_{j+1}^{(n)})T_j+d_{n-k-1}^{(n)}T_{n-k}+d_{n-k}^{(n)}T_{n-k+1}\big]
\cdot \big[\sum_{j=0}^kb_jT_j\big] \nonumber \\
&& + (c_1^{(n)}-\frac{c_{k-1}^{(n)}}{b_k}b_0)T_0+
\sum_{j=1}^{k-2}(c_{j-1}^{(n)}+c_{j+1}^{(n)}-\frac{c_{k-1}^{(n)}}{b_k}b_j)T_j
+(c_{k-2}^{(n)}-\frac{c_{k-1}^{(n)}}{b_k}b_{k-1})T_{k-1}, \nonumber
\end{eqnarray}
and subtract $T_{n-1}$ for identification of the $d_j^{(n+1)}$ and $c_j^{(n+1)}$,
\begin{equation}
T_{n+1}=2T_1T_n-T_{n-1}
=(\sum_{j=0}^{n-k+1}d_j^{(n+1)}T_j)(\sum_{j=0}^{k}b_jT_j)+
\mathop{{\sum}'}_{j=0}^{k-1}c_j^{(n+1)}T_j.
\end{equation}
\end{proof}

\begin{example}
For \eqref{eq.exa2}, we obviously have
\begin{equation}
c_0^{(1)}=
c_2^{(1)}=
c_0^{(2)}=
c_1^{(2)}=0,\quad
c_1^{(1)}=
c_2^{(2)}=1.
\end{equation}
in \eqref{eq.Tndiv}. The formulas \eqref{eq.dcRec1st}--\eqref{eq.dcReclast}
predict at $n=2$
\begin{equation}
d_0^{(3)}=\frac{1}{1/4},\quad
\frac{c_0^{(3)}}{2}=-\frac{78\frac{1}{2}\cdot 1}{1/4},\quad
c_1^{(3)}=1-\frac{-23\frac{1}{4}\cdot 1}{1/4}-1,\quad
c_2^{(3)}=-\frac{-1\frac{1}{2}\cdot 1}{1/4}.
\end{equation}
With these, \eqref{eq.anredu} gives at $n=3$
\begin{equation}
a_3=8+ (-314)\cdot a_0+93\cdot a_1+6\cdot a_2
\end{equation}
which is correct since
\begin{equation}
a_0\approx 0.02671606,
a_1\approx 0.00412578,
a_2\approx 0.00087916,
a_3\approx 0.00013030.
\label{eq.aexact03}
\end{equation}
The next step of the recursion is
\begin{equation}
d_0^{(4)}=\frac{6}{1/4},\quad
\frac{c_0^{(4)}}{2}=93-\frac{78\frac{1}{2}\cdot 6}{1/4},\quad
c_1^{(4)}=2\cdot(-314)+6-\frac{-23\frac{1}{4}\cdot 6}{1/4},\quad
c_2^{(4)}=93-\frac{-1\frac{1}{2}\cdot 6}{1/4}-1,
\end{equation}
\begin{equation}
a_4=48+ (-1791)\cdot a_0+(-64)\cdot a_1+128\cdot a_2
\end{equation}
which is also correct with
\begin{equation}
a_4\approx 0.00002159.
\label{eq.aexact4}
\end{equation}
\end{example}

\section{Approximation by the Truncated Chebyshev Series}\label{trunc.sec}
Approximations $\hat a_n$ to the $a_n$ of \eqref{eq.anDef}
may be calculated assuming that the $a_n$ are negligible beyond some index $N$:
\begin{equation}
\frac{1}{\sum_{j=0}^k b_jT_j(x)} \approx
\mathop{{\sum}'}_{n=0}^N \hat a_n T_n(x).
\end{equation}
If this equation is multiplied by $2\sum b_jT_j$, and we stay with
\eqref{eq.bbeyond} to keep the notation simple,
\begin{equation}
2\approx \mathop{{\sum}'}_{n=0}^N\hat a_n\sum_{l=0}^{k+n}
(b_{n-l}+b_{l+n}+b_{l-n})T_l.
\label{eq.ahat}
\end{equation}
If  the coefficients in front of $T_0$ to $T_N$ are set equal on
both sides, a system of linear equations for the $\hat a_n$ ensues:
\begin{equation}
\left(
\begin{array}{ccccc}
b_0 & b_1  & b_2 & b_3 & \ldots\\
b_1 & 2b_0+b_2 & b_1+b_3 & b_2+b_4 & \ldots \\
b_2 & b_1+b_3 & 2b_0+b_4 & b_1+b_5& \ldots  \\
b_3 & b_2+b_4 & b_1+b_5 & 2b_0+b_6& \ldots  \\
\vdots & \vdots & \vdots & \vdots & \ddots  \\
\end{array}
\right)
\cdot
\left(
\begin{array}{c}
\hat a_0 \\
\hat a_1 \\
\vdots \\
\vdots \\
\hat a_N \\
\end{array}
\right)
=
\left(
\begin{array}{c}
2 \\
0 \\
\vdots \\
\vdots \\
0
\end{array}
\right)
\label{eq.ahatMat}
\end{equation}
where the coefficient matrix $B_{r,c}$ is symmetric and has a band width of $2k+1$:
\begin{equation}
B_{r,c}=\left\{
\begin{array}{c@{\quad}c}
b_c, & r=0 \\
b_r, & c=0 \\
2b_0+b_{2r}, & r=c\neq 0 \\
b_{|r-c|}+b_{r+c}, & r\neq c,\quad c>0,\quad r>0 \\
\end{array}
\right.
\label{eq.BMat}
\end{equation}
This gives access to a set of approximate, low-indexed $a_n$ with no need
to evaluate integrals nor reference to the roots of $\sum b_jT_j$.

\begin{example}
Again for \eqref{eq.exa2}, the choice of $N=3$ yields the coefficient
vector
\begin{equation}
\hat a_0=0.02671602,
\hat a_1=0.00412567,
\hat a_2=0.00087845,
\hat a_3=0.00012696.
\end{equation}
At $N=4$, this improves to
\begin{equation}
\begin{array}{l}
\hat a_0=0.02671606,
\hat a_1=0.00412578,
\hat a_2=0.00087914,
\hat a_3=0.00013019,\\
\hat a_4=0.00002111,
\end{array}
\end{equation}
which is close to the exact results in \eqref{eq.aexact03}
and \eqref{eq.aexact4}.
At $N=5$, this improves further to
\begin{equation}
\begin{array}{l}
\hat a_0=0.02671606,
\hat a_1=0.00412578,
\hat a_2=0.00087916,
\hat a_3=0.00013029,\\
\hat a_4=0.00002158,\ldots
\end{array}
\end{equation}
\end{example}

\begin{remark}
The matrix in (\ref{eq.ahatMat}) is the approximate, square upper left
$(N+1)\times(N+1)$ submatrix of the ``exact'' solution. The approximate
solution obtained could be considered as if the terms
$\sum_{n=N+1}^{l+k}\hat a_n(b_{n-l}+b_{n+l}+b_{l-n})$ of (\ref{eq.ahat})
in the $l$'th row of
the system of linear equations had been neglected (as if the columns
$N+1$ up to $l+k$ had been chopped off).
The neglected sum is nonzero only if $l\ge N+1-k$.
An idea of an improvement of this algorithm is: reduce the $\hat a_n$
($N+1\le n \le l+k$) in the neglected terms via (\ref{eq.anredu})
to a linear combination of $\hat a_{1,\ldots,k}$, and re-introduce (add)
these components
(in)to the matrix---add the constant to the right hand side--- in these rows
$l\ge N+1-k$.
This update of the system of linear equations reduces the rank
of the matrix and does therefore not improve on what is obtained from (\ref{eq.ahatMat}).
\end{remark}

The algorithm may be extended to the division problem
of finding the $\hat a_n$ from given $f_n$ in
\begin{equation}
\frac{f(x)}{\sum_{j=0}^k b_jT_j(x)} \approx
\mathop{{\sum}'}_{n=0}^N \hat a_n T_n(x);\quad f(x)\equiv \mathop{{\sum}'}_{n=0}^\infty f_n T_n(x),
\label{eq.fndiv}
\end{equation}
with the right hand side in (\ref{eq.ahatMat}) replaced as follows:
\begin{equation}
\sum_{c=0}^N B_{r,c}\hat a_c = \left\{
\begin{array}{c@{,\quad}c}
f_r & r=0 ,\\
2f_r & r=1,2,3,\ldots \\
\end{array}
\right.
\label{eq.ahatN}
\end{equation}

\begin{example}\label{exam.sinxx}
The Chebyshev series of $f(x)=\sin(\frac{\pi}{2}x)/x$
starts with \cite{Clenshaw1954,Schonfelder}
\begin{eqnarray}
f_0&=&\pi\sum_{s=0}^\infty \frac{(-)^s}{(2s+1)(s!)^2}\left(\frac{\pi}{4}\right)^{2s}\approx 2.552557924804531760415274,\label{eq.fnsine}\\
f_2&\approx& -0.2852615691810360095702941, \\
f_4&\approx& 0.009118016006651802497767923, \\
f_6&\approx& -0.0001365875135419666724364765,\\
f_8&\approx& 0.000001184961857661690108290062,
\end{eqnarray}
\begin{equation}
f_n=\left\{
\begin{array}{c@{,\quad}c}
4(-)^{n/2}\sum_{s=n/2}^\infty J_{2s+1}(\pi/2) & n \quad \text{even}, \\
0 & n \quad \text{odd}.
\end{array}
\right.
\label{eq.sinxx}
\end{equation}
(Schonfelder \cite{Schonfelder} lists $2f_{2n}/\pi$ for $n\le 16$.)
If we approximate $f(x)$ by the polynomial $\mathop{{\sum}'}_{n=0}^4 f_n T_n(x)$,
calculation of the $\hat a_n$ in
\begin{equation}
\frac{f(x)}{\mathop{{\sum}'}_{j=0}^4 f_j T_j(x)}\approx \mathop{{\sum}'}_{n=0}^N \hat a_n T_n(x)
\end{equation}
via (\ref{eq.ahatN}) at $N=8$ predicts the relative error
\begin{eqnarray}
\mathop{{\sum}'}_{n=0}^N \hat a_n T_n(x)-1 
&\approx&
-6.74\cdot 10^{-8}T_0(x)-9.97\cdot 10^{-7}T_2(x) \nonumber \\
&& -1.23\cdot 10^{-5}T_4(x)
-1.09\cdot 10^{-4}T_6(x)-1.13\cdot 10^{-5}T_8(x).
\label{eq.sin4k}
\end{eqnarray}
\end{example}
\begin{remark}
The functional relations (\ref{eq.fndiv}) hold also for the shifted Chebyshev polynomials
$T^*(x)$:
\begin{equation}
\frac{f(x)}{\sum_{j=0}^k b_jT_j^*(x)} \approx
\mathop{{\sum}'}_{n=0}^N \hat a_n T_n^*(x);\quad f(x)\equiv \mathop{{\sum}'}_{n=0}^\infty f_n T_n^*(x),
\label{eq.fndivstar}
\end{equation}
\end{remark}

\section{Chebyshev Approximation for the Relative Error}\label{rel.sec}
The previous example of a truncated Chebyshev series had a maximum absolute
error estimated at $\sum_{n=6}^8 |f_n|\approx 0.000138$ if terms up
to $k=4$ were
retained, and the maximum relative error of the same polynomial was estimated at
$\mathop{{\sum}'}_{n=0}^8 |\hat a_n|-1\approx 0.000134$---dominated by
the $\hat a_6$ term
in (\ref{eq.sin4k}). To optimize the approximation of $f(x)$ for the
relative error in $[-1,1]$, one would rather like to find the $k+1$ coefficients
$b_j$ in (\ref{eq.fndiv}) which force the relative error to be close to
zero in the sense of
\begin{equation}
\hat a_0=2,\quad \hat a_1=\hat a_2=\hat a_3=\ldots=\hat a_k=0.
\label{eq.ak0}
\end{equation}
As an inversion of the problem of Sec.\ \ref{trunc.sec}, the matrix $B$ in (\ref{eq.ahatN})
is presumed unknown (up to some symmetry),
and the first $k+1$ elements of the vector $\hat a_c$ and all elements of $f_r$ are known.
The rationale is that removal of the ripples of $T_1(x)$ to $T_k(x)$ from the quotient
expansion leaves a quotient with an appropriate number of ``critical'' points
required by the alternating maximum theorem \cite{Cody,Nitsche,Veidinger}.
\begin{remark}
The case $r=0$ in (\ref{eq.ahatN}) in conjunction with (\ref{eq.ak0}) mandate
\begin{equation}
b_0=f_0/2 .
\label{eq.b0off0}
\end{equation}
\end{remark}
Finding the constituents $b_j$ of $B$ that solve
the bi-linear equation (\ref{eq.ahatN}) may
proceed with a vectorized first-order Newton method
as follows:
\begin{itemize}
\item
Chose a start solution $b_j$, for example the obvious
\begin{equation}
b_j=\left\{
\begin{array}{c@{\quad}c}
f_0/2, & j=0\\
f_j, & j=1,2,\ldots,k\\
\end{array}
\right.
\end{equation}
\item
Compute the $\hat a_n$ ($n=0,\ldots,N$) from $b_j$ by solving the
linear system of equations (\ref{eq.ahatN}).
\item
Compute an approximate $(N+1)\times k$ Jacobi matrix
\begin{equation}
J_{r,c}=
\left(
\begin{array}{ccccc}
\frac{\partial \hat a_0}{\partial b_1} & \frac{\partial \hat a_0}{\partial b_2}  & \ldots & \frac{\partial \hat a_0}{\partial b_k}\\
\frac{\partial \hat a_1}{\partial b_1} & \frac{\partial \hat a_1}{\partial b_2} & \ldots & \frac{\partial \hat a_1}{\partial b_k} \\
\vdots                                 & \vdots                                 & \ldots & \vdots                                 \\
\frac{\partial \hat a_N}{\partial b_1} & \frac{\partial \hat a_N}{\partial b_2} & \ldots & \frac{\partial \hat a_N}{\partial b_k} \\
\end{array}
\right)
\end{equation}
by partial derivation of the first $N+1$ equations
of  (\ref{eq.ahatN}) w.r.t.\ the $b_j$, i.e., by solving the
$k$ systems of $N+1$ linear equations
\begin{equation}
\sum_{c=0}^N B_{r,c}J_{c,j}=-
\left(
\begin{array}{cccccc}
\hat a_1 & \hat a_2 & \hat a_3 & \ldots & \hat a_{k-1} & \hat a_k\\
\hat a_0+\hat a_2 & \hat a_1+\hat a_3 & \hat a_2+\hat a_4 & \ldots & \ldots & \hat a_{k-1}+\hat a_{k+1}\\
\hat a_1+\hat a_3 & \hat a_0+\hat a_4 & \hat a_1+\hat a_5 & \ldots & \ldots & \hat a_{k-2}+\hat a_{k+2}\\
\vdots & \vdots & \vdots & \ldots & \vdots & \vdots \\
\hat a_{N-1} & \hat a_{N-2} & \hat a_{N-3} & \ldots & \hat a_{N-k-1} & \hat a_{N-k}
\end{array}
\right)
\end{equation}
for $r=0,\ldots,N$ and $j=0,\ldots,k-1$.
The column $\partial \hat a_c/\partial b_0$ of the Jacobi matrix is not calculated,
as $b_0$ is assumed fixed according to (\ref{eq.b0off0}).
\item
Compute the next iterated solution $b_j+\Delta_j$ ($j=1,2,\ldots,k$)
of the polynomial coefficients by solving the system of $k$ linear equations
\begin{equation}
\sum_{j=1}^k \frac{\partial \hat a_l}{\partial b_j} \Delta_j = -\hat a_l,\quad l=1,\ldots,k
\end{equation}
for the first-order differences $\Delta_j$. This equation is the first-order Taylor
expansion of $\hat a_l$ as a function of the $b_j$ set to the target (\ref{eq.ak0})
for this update.
The $k\times k$ coefficient matrix $\partial \hat a_l/\partial b_j$ is a
square submatrix of the Jacobi matrix calculated
in the previous step.
\item
Return to the second bullet for the next cycle until the $\hat a_0$ to
$\hat a_k$ are sufficiently close to (\ref{eq.ak0}).
\end{itemize}
\begin{remark}
This algorithm involves only $f_0$ to $f_N$, but no higher order
approximants to $f(x)$. It therefore adapts a polynomial of degree $k$ to
a polynomial of degree $N$.
\end{remark}
\begin{example}
The error terms (\ref{eq.sin4k}) for the polynomial $\sum_{j=0}^4b_jT_j(x)$
change to
\begin{eqnarray}
\mathop{{\sum}'}_{n=0}^N \hat a_n T_n(x)-1&\approx&
5.2\cdot 10^{-12}T_0(x)+4.7\cdot 10^{-11}T_2(x)+6.3\cdot 10^{-12}T_4(x) \nonumber \\
&& -1.08\cdot 10^{-4}T_6(x)-1.11\cdot 10^{-5}T_8(x)
\end{eqnarray}
after one Newton iteration, reducing the relative error to
$\mathop{{\sum}'}_{n=0}^8 |\hat a_n|-1\approx 0.000119$. During further iteration
cycles the relative error stays about the same because it is dominated
by $\hat a_6T_6(x)$ which is out of reach of the polynomial base with $k=4$.
\end{example}
\begin{example}
An IEEE ``single'' precision accuracy of $f(x)=\sin(\frac{\pi}{2}x)/x$ with
a relative error smaller than $2^{-24}\approx 6.0\cdot 10^{-8}$ needs
$k=8$.
Truncation of the Chebyshev series for $f(x)$ after $k=8$ yields
an estimated maximum absolute error of
$\sum_{n=k+1}^N |f_n|\approx 6.7\cdot 10^{-9}$ evaluated at $N=16$.
The relative error of the same polynomial is also
$\mathop{{\sum}'}_{n=0}^N |\hat a_n|-1\approx 6.7\cdot 10^{-9}$.
After four Newton iterations, this value drops to $5.9\cdot 10^{-9}$
with coefficients given in the following table---remaining very close
to those cited after (\ref{eq.fnsine}):

\begin{tabular}{cc}
$n$ & $b_n$ \\
\hline
0 & 1.276278962402265880207637 \\
2 & -0.2852615691810328617761446 \\
4 & 0.9118016006289075331306166$\cdot 10^{-2}$\\
6 & -0.1365874893444115901818408$\cdot 10^{-3}$\\
8 & 0.1184206224108742454613850$\cdot 10^{-5}$
\end{tabular}
\end{example}
\begin{example}\label{exam.cosxx}
$g(x)=\cos(\frac{\pi}{2}x)$ has the
expansion coefficients \cite{Clenshaw1954,Murnaghan,Schonfelder}
\begin{equation}
g_n=\left\{
\begin{array}{c@{,\quad}c}
2(-)^{n/2}J_n(\pi/2)& n\,\text{even}, \\
0& n\,\text{odd}.
\end{array}
\right.
\end{equation}

The approximation $g(x)\approx\mathop{{\sum}'}_{n=0}^k g_nT_n(x)$ has an
estimated maximum absolute error of 
$\sum_{n=k+1}^N |g_n|\approx 4.7\cdot 10^{-8}$ for the polynomial of degree
$k=8$ evaluated at $N=16$.
Because $g(x)$ is zero at both ends of the interval $[-1,1]$, the algorithm
does not find polynomials $\sum_{j=0}^kb_jT_j(x)$ with a uniformly convergent
Chebyshev expansion of the relative error---any $\hat a_n$ obtained depend
strongly on $N$. We therefore ``lift''
both zeros by looking at $f(x)=\cos(\frac{\pi}{2}x)/(1-x^2)$ instead, which has
the expansion coefficients
\begin{eqnarray}
f_0&=&\pi J_1(\pi/2), \\
f_2&=&f_0-2g_0, \\
f_n&=&2f_{n-2}-f_{n-4}-4g_{n-2},\quad n=4,6,8,\ldots,\\
f_n&=&0,\quad n\, \text{odd}.
\end{eqnarray}
Truncation of the Chebyshev series for $f(x)$ after $k=4$ yields
an estimated maximum absolute error of
$\sum_{n=k+1}^N |f_n|\approx 2.7\cdot 10^{-5}$ evaluated at $N=8$.
The relative error of the same polynomial is
$\mathop{{\sum}'}_{n=0}^N |\hat a_n|-1\approx 3.3\cdot 10^{-5}$.
After four Newton iterations, this value drops to $3.1\cdot 10^{-5}$
with coefficients given in the following table:

\begin{tabular}{cc}
$n$ & $b_n$ \\
\hline
0 & 0.8903651967922106931461297 \\
2 & -0.1072744347398521266520654 \\
4 & 0.002332103968386755210894198 \\
\end{tabular}
\end{example}
\begin{example}
The coefficients of the Chebyshev series of $\arcsin x$
and $(\arcsin x)/x$ (App.\ \ref{app.arcsin})
are slowly descending. The infinite slope of
$\arcsin x$ at $x=\pm 1$ renders both series inefficient, so we
turn to $\frac{1}{x}\arcsin\frac{x}{\sqrt{2}}$ instead as configured in
(\ref{eq.asin2f}).
Keeping terms up to $f_{36}$ yields an estimated maximum absolute error of
$\sum_{n=38}^N |f_n|\approx 1.4\cdot 10^{-17}$ evaluated at $N=108$.
The relative error of the same polynomial is
$\approx 1.9\cdot 10^{-17}$.
After four Newton iterations, this value drops only slightly to
$1.8\cdot 10^{-17}$; obviously, there is not much room to improve the
polynomial representation w.r.t.\ an optimized {\em relative} error in
cases where the amplitude of the function is small over the $x$-interval.
\end{example}

\begin{example}
The expansion for $\exp(x)$ in $-1\le x\le 1$ reads \cite[(9.6.19)]{AS}
\cite[(3.4.1e)]{Rivlin}\cite[p.\ 69]{Fox}\cite[(33)]{ElliottMathComp19}
\begin{equation}
\exp(x)= 2\mathop{{\sum}'}_{n=0}^\infty I_n(1)T_n(x).
\end{equation}
Truncation after $k=14$ yields an estimated maximum absolute error of
$\sum_{n=k+1}^N |f_n|\approx 4.9\cdot 10^{-17}$ evaluated at $N=42$.
The relative error of the same polynomial is
$\mathop{{\sum}'}_{n=0}^N |\hat a_n|-1\approx 1.3\cdot 10^{-16}$.
After four Newton iterations, this value drops to $7.5\cdot 10^{-17}$
with coefficients given in the following table,
also listed as $f(x)\approx \sum_{n=0}^N d_nx^n$:

\begin{tabular}{ccc}
$n$ & $b_n$ & $d_n$\\
\hline
0 & 1.2660658777520083355982446 & 1.00000000000000002107745526254 \\
1 & 1.1303182079849700544153921 & 1.00000000000000063548946139343 \\
2 & 0.2714953395340765623657051 & 0.499999999999997953936666685291 \\
3 & 0.4433684984866380495257150$\cdot 10^{-1}$ & 0.1666666666666422610320391 \\
4 & 0.5474240442093732650276168$\cdot 10^{-2}$ & 0.4166666666669875817272051$\cdot 10^{-1}$ \\
5 & 0.5429263119139437503621352$\cdot 10^{-3}$ & 0.8333333333602639662588442$\cdot 10^{-2}$ \\
6 & 0.4497732295429514665443872$\cdot 10^{-4}$ & 0.1388888888702869286166025$\cdot 10^{-2}$ \\
7 & 0.3198436462401990501334121$\cdot 10^{-5}$ & 0.1984126971086418099245159$\cdot 10^{-3}$ \\
8 & 0.1992124806672795001043316$\cdot 10^{-6}$ & 0.2480158780231612103680909$\cdot 10^{-4}$ \\
9 & 0.1103677172551632915777862$\cdot 10^{-7}$ & 0.2755735152373104259316644$\cdot 10^{-5}$ \\
10 & 0.5505896079551881657982078$\cdot 10^{-9}$ & 0.2755725369287090362239172$\cdot 10^{-6}$ \\
11 & 0.2497956604792065959497342$\cdot 10^{-10}$ & 0.2504783672757589754944252$\cdot 10^{-7}$ \\
12 & 0.1039151254481832513826561$\cdot 10^{-11}$ & 0.2088034159586738951818317$\cdot 10^{-8}$ \\
13 & 0.3990676874210170341122722$\cdot 10^{-13}$ & 0.1634581247676485771723867$\cdot 10^{-9}$ \\
14 & 0.1400237499722866786358850$\cdot 10^{-14}$ & 0.1147074559772972471385170$\cdot 10^{-10}$
\end{tabular}
\end{example}

\begin{example}
The expansion for $J_0(\frac{\pi}{2}x)$ in $-1\le x\le 1$ reads
\cite[6.681.5]{GR}
\begin{equation}
J_0(\frac{\pi}{2}x)= 2\mathop{{\sum}'}_{n=0,2,4,\ldots}^\infty
(-)^{n/2}J_{n/2}^2(\frac{\pi}{4})T_n(x).
\end{equation}
Truncation after $T_{16}(x)$ yields an estimated maximum absolute error of
$\sum_{n=k+1}^N |f_n|\approx 7.3\cdot 10^{-19}$ evaluated at $N=48$.
The relative error of the same polynomial is
$\mathop{{\sum}'}_{n=0}^N |\hat a_n|-1\approx 1.6\cdot 10^{-18}$.
After four Newton iterations, this value drops to $1.3\cdot 10^{-18}$
with coefficients given in the following table:

\begin{tabular}{ccc}
$n$ & $b_n$ & $d_n$\\
\hline
0 & 0.7252769164405135618043045 & 0.9999999999999999991311745 \\
2 & -0.2638108118461404734713153 & -0.6168502750680847778603892 \\
4 & 0.1072184541022420669256084$\cdot 10^{-1}$ & 0.9512606546288948620024320$\cdot 10^{-1}$ \\
6 & -0.1885687642135952967199171$\cdot 10^{-3}$ & -0.6519837738512518004083602$\cdot 10^{-2}$ \\
8 & 0.1845983728936489887451460$\cdot 10^{-5}$ & 0.2513602312234872245916252$\cdot 10^{-3}$ \\
10 & -0.1150537142155094251800350$\cdot 10^{-7}$ & -0.6202064609606906421245435$\cdot 10^{-5}$ \\
12 & 0.4965029850154789447530764$\cdot 10^{-10}$ & 0.1062698637612363296679714$\cdot 10^{-6}$ \\
14 & -0.1571252252452718608949964$\cdot 10^{-11}$ & -0.1336990135568532922581048$\cdot 10^{-8}$ \\
16 & 0.3800986508122698831881511$\cdot 10^{-15}$ & 0.1245507258981645953230933$\cdot 10^{-10}$
\end{tabular}
\end{example}

A set of $b_j$ in 
\begin{equation}
R(x)\equiv \frac{f(x)}{\sum_{j=0}^k b_jT_j(x)}-1
\label{eq.R}
\end{equation}
found that way is also a starting point to calculate
the solution with the minimax property of the relative error: This locates
the local minima and maxima of $R(x)$, computes the mean of their absolute
values, and iteratively adjusts the $b_j$ such that the absolute values of the
new alternating extrema equal that mean. The corrections
$\Delta_j$ to the $b_j$ can be computed by 
expansion of (\ref{eq.R}) to first order in $\Delta_j$
keeping the abscissa of the extrema fixed, which
ends up in a linear system of equations for the $\Delta_j$.

\begin{example}\label{exam.sinxxD}
An IEEE ``double'' precision accuracy of $f(x)=\sin(\frac{\pi}{2}x)/x$ with
a relative error smaller than $2^{-53}\approx 1.1\cdot 10^{-16}$ needs
$k=16$.
Truncation of the Chebyshev series of Example (\ref{exam.sinxx}) for
$f(x)$ after $k=16$ yields an estimated maximum absolute error of
$\sum_{n=k+1}^N |f_n|\approx 4.1\cdot 10^{-19}$ evaluated at $N=32$.
The relative error of the same polynomial is
$\mathop{{\sum}'}_{n=0}^N |\hat a_n|-1\approx 3.8\cdot 10^{-19}$.
After four Newton iterations, this value drops to $3.5\cdot 10^{-19}$
with coefficients $b_n$ given in the following table:

\begin{tabular}{ccc}
$n$ & $b_n$ & $d_n$\\
\hline
0 & 1.276278962402265880207637 & 1.5707963267948966188688195 \\
2 & -0.2852615691810360095702941 & -0.6459640975062461962319336 \\
4 & 0.9118016006651802497767923$\cdot 10^{-2}$ & 0.7969262624616554097627533$\cdot 10^{-1}$ \\
6 & -0.1365875135419666724364765$\cdot 10^{-3}$ & -0.4681754135303468240882506$\cdot 10^{-2}$ \\
8 &   0.1184961857661690108288872$\cdot 10^{-5}$ & 0.1604411847100114088031881$\cdot 10^{-3}$ \\
10 & -0.6702791603827441081706121$\cdot 10^{-8}$ & -0.3598843013917326159520456$\cdot 10^{-5}$ \\
12 &  0.2667278599017903283863443$\cdot 10^{-10}$ & 0.5692135656122429901944357$\cdot 10^{-7}$ \\
14 & -0.7872922004615709018594325$\cdot 10^{-13}$ & -0.6684369436484103757933363$\cdot 10^{-9}$ \\
16 &  0.1791929094718284072119916$\cdot 10^{-15}$ & 0.5871793257572873247522307$\cdot 10^{-11}$
\end{tabular}
The actual relative error of this approximation is shown in Fig.\ \ref{err.ps} as a continuous
line, with a maximum of $2.9\cdot 10^{-19}$. The dashed line with a relative error of
$2.6\cdot 10^{-19}$ in comparison results from 
further minimax optimization with coefficients shown in the next table:
\begin{figure}[h]
\includegraphics[scale=0.5]{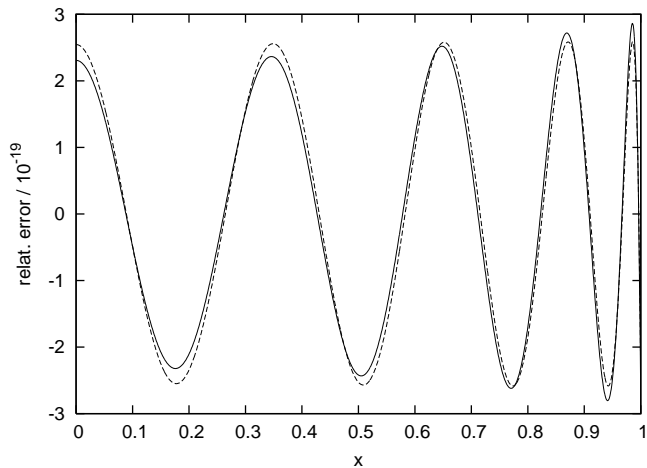}
\caption{The relative error $R(x)$ to $f(x)=\sin(\pi x/2)/x$ for both tabulated
parameter sets of $\sum_{n=0,2,\ldots}^{16} b_n T_n(x)$ of Example \ref{exam.sinxxD}.
}
\label{err.ps}
\end{figure}

\begin{tabular}{ccc}
$n$ & $b_n$ & $d_n$ \\
\hline
0 & 1.2762789624022658802075437 & 1.5707963267948966188314659 \\
2 & -0.2852615691810360095705230 & - 0.6459640975062461915471363 \\
4 & 0.9118016006651802497528156$\cdot 10^{-2}$ &  0.7969262624616544421893744$\cdot 10^{-1}$ \\
6 & -0.1365875135419666726405733$\cdot 10^{-3}$ & -0.4681754135302704719117724$\cdot 10^{-2}$\\
8 &   0.1184961857661689920542732$\cdot 10^{-5}$ &  0.1604411847070460989830944$\cdot 10^{-3}$\\
10 & -0.6702791603827612608171959$\cdot 10^{-8}$ & -0.3598843007652658555526627$\cdot 10^{-5}$\\
12 &  0.2667278599019855592489579$\cdot 10^{-10}$ &  0.5692134921914455833455723$\cdot 10^{-7}$\\
14 & -0.7872921659616258733890169$\cdot 10^{-13}$ & -0.6684324580312975131354658 $\cdot 10^{-9}$\\
16 &  0.1791589025538146793760922$\cdot 10^{-15}$ & 0.5870678918883399413795788$\cdot 10^{-11}$
\end{tabular}
\end{example}

\begin{example}
As an example for (\ref{eq.fndivstar}), consider
$\exp(x)=\mathop{{\sum}'}_{n=0}^\infty f_nT^*_n(x)$
over $0\le x\le 1$
\cite[(4.2.48)]{AS}\cite{Clenshaw1954,KhajahCMA27}. The $f_n$ are represented
via \cite[(9.6.26]{AS} through modified Bessel Functions $I_n$,
\begin{equation}
f_n=2\sqrt{e}I_n(1/2);\quad f_{n+1}=-4nf_n+f_{n-1},
\end{equation}

\begin{tabular}{cc}
$n$ & $f_n$ \\
\hline
0 & 3.506775308754180791443893 \\
1 & 0.8503916537808109665352350 \\
2 & 0.1052086936309369253029528 \\
3 & 0.008722104733315564111612874 \\
4 & 0.0005434368311501559635982758 \\
5 & 0.00002711543491306869404046064
\end{tabular}

Truncation of the Chebyshev series for $f(x)$ after $k=3$ yields
an estimated maximum absolute error of
$\sum_{n=k+1}^N |f_n|\approx 5.7\cdot 10^{-4}$ evaluated at $N=9$.
The relative error of the same polynomial is
$\mathop{{\sum}'}_{n=0}^N |\hat a_n|-1\approx 5.1\cdot 10^{-4}$.
After four Newton iterations, this value drops to $4.0\cdot 10^{-4}$
with coefficients given in the following table:

\begin{tabular}{cc}
$n$ & $b_n$ \\
\hline
0 & 1.753387654377090395721946 \\
1 & 0.8503902561425088936327743 \\
2 & 0.1051918520893768747555014 \\
3 & 0.008587089960927766771654559
\end{tabular}

If we proceed to $k=12$ at $N=36$, the estimated maximum relative error
becomes $6.1\cdot 10^{-18}$ with the following coefficients:

\begin{tabular}{ccc}
$n$ & $b_n$ & $d_n$ \\
\hline
0 & 1.7533876543770903957219464 & 1.0000000000000000060373678 \\
1 & 0.8503916537808109665352350 & 0.9999999999999978889799411\\
2 & 0.1052086936309369253029528 & 0.5000000000001216148194572 \\
3 & 0.8722104733315564111612874$\cdot 10^{-2}$ & 0.1666666666639271874501180 \\
4 & 0.5434368311501559635982758$\cdot 10^{-3}$ & 0.4166666669859109153386033$\cdot 10^{-1}$ \\
5 & 0.2711543491306869404045765$\cdot 10^{-4}$ & 0.8333333112815145481691497$\cdot 10^{-2}$ \\
6 & 0.1128132888782082788967416$\cdot 10^{-5}$ & 0.1388889862738933258163839$\cdot 10^{-2}$ \\
7 & 0.4024558229870710027066467$\cdot 10^{-7}$ & 0.1984098287973665146421103$\cdot 10^{-3}$ \\
8 & 0.1256584418283842256517024$\cdot 10^{-8}$ & 0.2480734627092463176804164$\cdot 10^{-4}$ \\
9 & 0.3488091362080888722258141$\cdot 10^{-10}$ & 0.2747848541489261879291146$\cdot 10^{-5}$ \\
10 & 0.8715278679388174731063544$\cdot 10^{-12}$ & 0.2827881515524984459349078$\cdot 10^{-6}$ \\
11 & 0.1979783472020383084286900$\cdot 10^{-13}$ & 0.2086709669366350082217004$\cdot 10^{-7}$ \\
12 & 0.4103178180353125619414324$\cdot 10^{-15}$ & 0.3441995330913567239602395$\cdot 10^{-8}$
\end{tabular}
Equilibration of the local extrema with the following coefficients
reduces this error to $5.0\cdot 10^{-18}$:

\begin{tabular}{ccc}
$n$ & $b_n$ & $d_n$ \\
\hline
0 & 1.7533876543770903961757996 & 1.0000000000000000049913878 \\
1 & 0.8503916537808109674449984 & 0.9999999999999982006556032\\
2 & 0.1052086936309369262175803 & 0.5000000000001063630793784 \\
3 & 0.8722104733315565035195027$\cdot 10^{-2}$ & 0.1666666666642173677701902 \\
4 & 0.5434368311501568988037582$\cdot 10^{-3}$ & 0.4166666669575619866774579$\cdot 10^{-1}$ \\
5 & 0.2711543491306964588901369$\cdot 10^{-4}$ & 0.8333333129071360606767972$\cdot 10^{-2}$ \\
6 & 0.1128132888783054546376918$\cdot 10^{-5}$ & 0.1388889803905621871712292$\cdot 10^{-2}$ \\
7 & 0.4024558229970401854905218$\cdot 10^{-7}$ & 0.1984099684263107292542263$\cdot 10^{-3}$ \\
8 & 0.1256584419307581158766302$\cdot 10^{-8}$ & 0.2480712594971168345247335$\cdot 10^{-4}$ \\
9 & 0.3488091466142981067488685$\cdot 10^{-10}$ & 0.2748077480706561519632930$\cdot 10^{-5}$ \\
10 & 0.8715288225355426665019433$\cdot 10^{-12}$ & 0.2826376902452534679601836$\cdot 10^{-6}$ \\
11 & 0.1979820112783685973416909$\cdot 10^{-13}$ & 0.2092376435267840500466468$\cdot 10^{-7}$ \\
12 & 0.4092071997914099904014169$\cdot 10^{-15}$ & 0.3432678789827820176761249$\cdot 10^{-8}$
\end{tabular}
\end{example}

\section{Summary}
Besides some generic algorithms to compute the Chebyshev series
of inverse polynomials, there are two specific aspects that
facilitate this task: (i) the expansion coefficients can be
derived from the partial fractions of the inverse polynomial.
(ii) Expansion coefficients with indices larger than the
polynomial degree are recursively linked to those
of lower order. (iii) An algorithm has been presented which
derives a polynomial of a given degree such that the first
terms of the Chebyshev expansion of the relative error
of a given function represented by this polynomial vanish.

\appendix
\section{Chebyshev Series of $\ln(1+x)$}
The integral representation
\begin{equation}
\ln(1+x)=\int\frac{dx}{1+x}
\end{equation}
and term-by-term integration of (\ref{eq.T1x}) on the r.h.s.\ with
(\ref{eq.Tstarintx}) yield that the Chebyshev coefficients of
\begin{equation}
f(x)=\ln(1+x)\equiv
\mathop{{\sum}'}_{n=0}^\infty f_nT_n^*(x)\,\quad 0\le x\le 1,
\label{eq.Tlogx}
\end{equation}
obey
\begin{equation}
2nf_n=a_{n+1,1}(-3)-a_{n-1,1}(-3),\quad n\ge 1,
\end{equation}
explicitly \cite[p.\ 88]{Fox}\cite{FraserJACM12}
\begin{equation}
f_n=\frac{2(-)^{n+1}}{n(3+2\sqrt{2})^n},\quad n\ge 1
\end{equation}
from (\ref{eq.an-3}), as tabulated in \cite[4.1.45]{AS}. The missing $f_0$
is
\begin{equation}
f_0=\frac{2}{\pi}\int_0^1 \frac{\ln(1+x)}{\sqrt{x(1-x)}}dx=2\ln\frac{3+2\sqrt{2}}{4},
\end{equation}
because insertion of $x=1$ in (\ref{eq.Tlogx}) yields
\begin{equation}
f_0=2\left(f(1)-\sum_{n=1}^\infty f_n\right)
\end{equation}
and $f(1)=\ln 2$ and $\sum_{n=1}^\infty f_n=2\ln(1+\frac{1}{3+2\sqrt2})$
via \cite[4.1.24]{AS}.

\section{Chebyshev Series of $\arctan x$}
Integrating (\ref{eq.1xx}) over $x$ with
\begin{equation}
\int \frac{1}{1+x^2}dx=\arctan x
\end{equation}
and (\ref{eq.Tintx}) we get \cite[p.\ 89]{Fox}\cite{FraserJACM12}
\begin{equation}
\arctan x=2\sum_{j=1,3,5,7,\ldots}\frac{(-)^{\lfloor j/2\rfloor}}{j(1+\sqrt{2})^j}T_j(x),
\quad -1\le x\le 1,
\label{eq.atanT}
\end{equation}
in particular at $x=1$
\begin{equation}
\frac{\pi}{8}=\sum_{j=1,3,5,7,\ldots}\frac{(-)^{\lfloor j/2\rfloor}}{j(1+\sqrt{2})^j}.
\end{equation}
From (\ref{eq.atanT}) and (\ref{eq.fxx}), the coefficients of
\begin{equation}
\frac{\arctan x}{x}\equiv \mathop{{\sum}'}_{n=0,2,4,6,\ldots}g_nT_n(x)
\end{equation}
as listed in \cite[(4.4.50)]{AS}\cite{Clenshaw1954}
follow recursively, where $g_0=2\ln(1+\sqrt{2})$
is obtained via \cite[4.531.12]{GR}.

\section{Chebyshev Series of $\arcsin x$}\label{app.arcsin}
The series of $\arcsin x=\sum_{n=1,3,5,\ldots}^\infty g_nT_n(x)$
starts with $g_1=4/\pi$. A combination of \cite[(4.4.58)]{AS}, 
(\ref{eq.Tnmsum}), (\ref{eq.Tintx}) and \cite[(3.4.1d)]{Rivlin}
\begin{equation}
\sqrt{1-x^2}=\frac{4}{\pi} \mathop{{\sum}'}_{
\genfrac{}{}{0pt}{}{n=0}{n\, \text{even}}
}^\infty \frac{1}{1-n^2}T_n(x)
\end{equation}
yields
\begin{equation}
(n+1)g_{n+1}+(n-1)g_{n-1}=\frac{8}{\pi}\frac{n}{n^2-1},\quad n=2,4,6,.\dots
\end{equation}
in this case, which can be unwound as $g_n=4/(\pi n^2)$.
To find a formulation with controlled {\em relative} error,
we would switch to $h(x)=(\arcsin x)/x= \mathop{{\sum}'}_{n=0}^\infty h_nT_n(x)$
to remove the zero in the spirit of example \ref{exam.cosxx}.
With (\ref{eq.fxx}), the expansion coefficients are
\begin{eqnarray}
h_0&=&8\beta(2)/\pi, \\
h_{n+2}&=& -h_n+\frac{8}{\pi(n+1)^2},
\end{eqnarray}
where $\beta(2)\approx 0.915965594177219015054603515$ is Catalan's constant
\cite[Tab 23.3]{AS}\cite[0.234.3]{GR}.

\section{Chebyshev Series of $\arcsin(x/\surd 2)$}\label{app.arcsin2}
The coefficients of
\begin{equation}
\arcsin(x/\sqrt{2})=\sum_{n=1,3,5,\ldots}^\infty k_nT_n(x)
\end{equation}
are found by partial integration
of $\int \arcsin([\cos\theta]/\surd 2)\cos(n\theta)d\theta$
\begin{equation}
k_n
=\frac{1}{n\pi}\left[
\int_0^\pi \frac{\cos[(n-1)\theta]}{\sqrt{2-\cos^2 \theta}}d\theta
-\int_0^\pi \frac{\cos[(n+1)\theta]}{\sqrt{2-\cos^2 \theta}}d\theta
\right],\quad n\, \text{odd}
\label{eq.acosf}
\end{equation}
where
\begin{equation}
\int_0^{\pi}
\frac{\cos(2m\theta)}{\sqrt{2-\cos^2\theta}}d\theta
=
G_{2m}
=
\left\{
\begin{array}{ll}
\sqrt{2}F(\frac{1}{\sqrt{2}})\approx 2.6220575542921198104648395899, & m=0, \\
\sqrt{2}[3F(\frac{1}{\sqrt{2}})-4E(\frac{1}{\sqrt{2}})]\approx 0.22577708482093539558499460534, & m=1, \\
\end{array}
\right.
\label{eq.G02}
\end{equation}
are Complete Elliptic Integrals.
To find a recurrence for these
\begin{equation}
G_s\equiv \int_{-1}^1 \frac{T_{s}(x)}{\sqrt{(2-x^2)(1-x^2)}}dx,
\end{equation}
we apply the method of \cite[(17.1.4)]{AS} to the quartic
$y^2\equiv (2-x^2)(1-x^2)$, with
$d(y T_s(x))/dx=y (dT_s(x)/dx)+T_s \frac{1}{2y}(T_3(x)-3T_1(x))$,
insert (\ref{eq.Tderi}) for the derivative on the r.h.s,
replace the first $y$ on the r.h.s.\ by
$y^2/y=(T_4/8-T_2+7/8)/y$, expand all products with (\ref{eq.Tnmsum}), and
finally insert the upper limit $x=1$ where $y(x)T_s(x)=0$:
\begin{equation}
G_{s+3}+G_{|s-3|}-3(G_{s+1}+G_{|s-1|})
+ \frac{s}{2}\mathop{{\sum}'}_{\genfrac{}{}{0pt}{}{l=0}{l-s\, \text{odd}}}^{s-1}
\left[
G_{l+4}+G_{|l-4|} -8(G_{l+2}+G_{|l-2|}) +14G_l
\right]
=0 .
\end{equation}
Inserting $s=1,3,5$ and 7, for example, yields
\begin{eqnarray}
3G_4-12G_2+G_0&=&0, \\
5G_6-27G_4+15G_2-G_0&=&0, \\
7G_8-41G_6+29G_4-3G_2&=&0, \\
9G_{10}-55G_8+43G_6-5G_4&=&0,
\end{eqnarray}
and generates $k_1$ to $k_{9}$ in (\ref{eq.acosf}) from $G_0$ and
$G_2$
shown in (\ref{eq.G02}).
A slowly converging series expansion is also known \cite[806.01]{Byrd}.
With (\ref{eq.fxx}) we find the
coefficients $f_n=2k_{n-1}-f_{n-2}$ for
\begin{equation}
\frac{1}{x}\arcsin\frac{x}{\sqrt{2}}= \mathop{{\sum}'}_{n=0,2,4,\ldots}^\infty f_nT_n(x),\quad -1\le x\le 1,
\label{eq.asin2f}
\end{equation}
starting at
\begin{eqnarray}
f_0&=&\frac{4}{\pi}\int_0^1 \frac{\arcsin\frac{x}{\sqrt{2}}}{x\sqrt{1-x^2}}dx
=\sqrt{2}\sum_{l=0}^\infty \frac{[(2l-1)!!]^2}{2^l(2l+1)[(2l)!!]^2} \nonumber\\
&=& \frac{\pi}{2}-\frac{4}{\pi}
\sum_{q=1}^\infty \frac{1}{4q-1}\left[G_{4q-2}-G_{4q}\right]
\approx
1.4866664932871034689603296833 .
\end{eqnarray}
The four coefficients $\alpha_i$ that span
\begin{equation}
k_{n-1}=\frac{\surd 2}{\pi}\left[
\alpha_1 K(\frac{1}{\sqrt{2}})
+\alpha_2 E(\frac{1}{\sqrt{2}})
\right], \quad
f_{n}=\frac{\surd 2}{\pi}\left[
\alpha_3 K(\frac{1}{\sqrt{2}})
+\alpha_4 E(\frac{1}{\sqrt{2}})
\right]+(-)^{[n/2]}f_0,
\end{equation}
start as follows:

\begin{tabular}{c|cc|cc}
$n$ & $\alpha_1$ & $\alpha_2$ & $\alpha_3$ & $\alpha_4$  \\
\hline
2 & -2 & 4 & -4 & 8 \\
4 & -26/9 & 4 & -16/9 & 0 \\
6 & -638/75 & 292/25 & -3428/225 & 584/25 \\
8 & -22702/735 & 212/5 & -513088/11025 & 1536/25 \\
10 & -23722/189 & 4652/27 & -6763436/33075 & 191128/675 \\
12 & -463174/847 & 2252/3 & -3558618544/4002075 & 822272/675 \\
14 & -162508858/65065 & 8691484/2535 & -2777152623884/676350675 & 643269592/114075 \\
\end{tabular}

Because $T_{2j}(x)=T_j^*(x^2)$, the following numbers coincide with
\cite[(4.4.51)]{AS} up to a factor $\surd 2$:

\begin{tabular}{c|c|c|c}
$n$ & $f_n$ & $n$ & $f_n$ \\
\hline
2 & 0.3885303371652290716432228$\cdot 10^{-1}$ \\
4 & 0.2885441422084471126676825$\cdot 10^{-2}$ &
6 & 0.2884218334475536563483289$\cdot 10^{-3}$ \\
8 & 0.3322367192785279209254231$\cdot 10^{-4}$ &
10 & 0.4158477878052832866177270$\cdot 10^{-5}$ \\
12 & 0.5496504525974164467345493$\cdot 10^{-6}$ &
14 & 0.7550078449371525934251585$\cdot 10^{-7}$ \\
16 & 0.1067193805629843129424091$\cdot 10^{-7}$ &
18 & 0.1542180379281470021561106$\cdot 10^{-8}$ \\
20 & 0.2268114598545151963877153$\cdot 10^{-9}$ &
22 & 0.3383885639342775871004709$\cdot 10^{-10}$ \\
24 & 0.5108937524377197224216916$\cdot 10^{-11}$ &
26 & 0.7791139213632464421446539$\cdot 10^{-12}$ \\
28 & 0.1198378589352895337866326$\cdot 10^{-12}$ &
30 & 0.1856972621821342234640637$\cdot 10^{-13}$ \\
32 & 0.2896189154386304361020997$\cdot 10^{-14}$ &
34 & 0.4542792886328823081478511$\cdot 10^{-15}$ \\
36 & 0.7161678029265506176831289$\cdot 10^{-16}$ &
38 & 0.1134144256904559996509711$\cdot 10^{-16}$
\end{tabular}

\section{Chebyshev Series of $\psi(x+2)$}\label{app.psi}

An expansion of the Digamma function \cite{SpougeSIAMJ31} is \cite[(6.3.16)]{AS}
\begin{equation}
\psi(2+x)=1-\gamma+x\sum_{k=2}^\infty\frac{1}{k(x+k)},
\end{equation}
where $\gamma\approx 0.5772$ is Euler's constant. Employing $a_{n,1}(-k)$
of (\ref{eq.Haseg}),
\begin{equation}
\frac{1}{x+k}=-\frac{1}{-k-x}=\frac{2}{\sqrt{k^2-1}}
\mathop{{\sum}'}_{n=0}^\infty \frac{(-)^n}{(k+\sqrt{k^2-1})^n}T_n(x),
\quad -1\le x\le 1,
\end{equation}
\begin{equation}
\psi(2+x)=1-\gamma+2x\sum_{k=2}^\infty
\frac{1} {k\sqrt{k^2-1}}
\mathop{{\sum}'}_{n=0}^\infty \frac{(-)^n}{(k+\sqrt{k^2-1})^n}T_n(x).
\label{eq.psiT}
\end{equation}
The auxiliary definition
\begin{equation}
K_n\equiv \sum_{k=2}^\infty \frac{1}{k\sqrt{k^2-1}(k+\sqrt{k^2-1})^n},
\quad n=0,1,2,\ldots
\label{eq.Kdef}
\end{equation}
turns (\ref{eq.psiT}) with the aid of (\ref{eq.Tnmsum}) into
\begin{equation}
\psi(x+2)=(1-\gamma-K_1)T_0(x)-\sum_{n=1}^\infty (-)^n(K_{n-1}+K_{n+1})T_n(x),\quad -1\le x\le 1,
\end{equation}
where
\begin{equation}
K_n+K_{n+2}= 2\sum_{k=2}^\infty \frac{1}{\sqrt{k^2-1}(k+\sqrt{k^2-1})^{n+1}},
\quad n=0,1,2,\ldots
\label{eq.Kdiff}
\end{equation}
Alternatives to the slowly converging original series (\ref{eq.Kdef}) at small $n$ are
obtained in terms of the Riemann Zeta function $\zeta$ after
reducing the fraction in (\ref{eq.Kdef}) and/or (\ref{eq.Kdiff})
by $k-\sqrt{k^2-1}$,
\begin{eqnarray}
K_0 &=&
\sum_{k=2}^\infty \frac{1}{k^2}\left(1-\frac{1}{k^2}\right)^{-1/2}
=
\sum_{l=0}^\infty (-)^l \genfrac{(}{)}{0pt}{}{-1/2}{l}[\zeta(2l+2)-1], \\
K_1 &=&
\sum_{k=2}^\infty \frac{k-\sqrt{k^2-1}}{k\sqrt{k^2-1}}
=
\sum_{l=1}^\infty (-)^l \genfrac{(}{)}{0pt}{}{-1/2}{l}[\zeta(2l+1)-1],\\
K_0+K_2 &=&
2\sum_{l=1}^\infty (-)^l \genfrac{(}{)}{0pt}{}{-1/2}{l}[\zeta(2l)-1], \\
K_1+K_3 &=&
2\sum_{l=2}^\infty (-)^l \left( \genfrac{(}{)}{0pt}{}{-1/2}{l}
+ \genfrac{(}{)}{0pt}{}{1/2}{l} \right) [\zeta(2l-1)-1], \\
K_n+K_{n+2} &=&
2\sum_{l=[(n+3)/2]}^\infty (-)^l
[\zeta(2l-n)-1]
\sum_{s=0}^{n+1} \genfrac{(}{)}{0pt}{}{n+1}{s}
 \genfrac{(}{)}{0pt}{}{(s-1)/2}{l} .
\end{eqnarray}

\begin{tabular}{c|c|c|c}
$n$ & $K_n$ & $n$ & $K_n$ \\
\hline
0 & 0.6942240199692270653811973 &
1 & 0.1181923495113155830503315 \\
2 & 0.2615575442260127035429158$\cdot 10^{-1}$ &
3 & 0.6357242927298094244957032$\cdot 10^{-2}$ \\
4 & 0.1613702909326556648518537$\cdot 10^{-2}$ &
5 & 0.4189942166841513997803225$\cdot 10^{-3}$ \\
6 & 0.1101726048982138724638504$\cdot 10^{-3}$ &
7 & 0.2918277395837793537278094$\cdot 10^{-4}$ \\
8 & 0.7763995103341854698876680$\cdot 10^{-5}$ &
9 & 0.2071120322602199079344235$\cdot 10^{-5}$ \\
10 & 0.5534045978754736410165904$\cdot 10^{-6}$ &
11 & 0.1480224417758054637706871$\cdot 10^{-6}$ \\
12 & 0.3961806941781982189370558$\cdot 10^{-7}$ &
13 & 0.1060807013890109056491206$\cdot 10^{-7}$ \\
14 & 0.2841134565373781348071928$\cdot 10^{-8}$ &
15 & 0.7610594780500236739360477$\cdot 10^{-9}$ \\
16 & 0.2038876075855356359426642$\cdot 10^{-9}$ &
17 & 0.5462507275750785409310247$\cdot 10^{-10}$ \\
18 & 0.1463563991421891531670298$\cdot 10^{-10}$ &
19 & 0.3921418684181661587649434$\cdot 10^{-11}$ \\
20 & 0.1050708536652289553110610$\cdot 10^{-11}$ &
21 & 0.2815309431326615497301838$\cdot 10^{-12}$ \\
22 & 0.7543503527420115661214215$\cdot 10^{-13}$ &
23 & 0.2021259323594124792356794$\cdot 10^{-13}$ \\
24 & 0.5415919982188370035264436$\cdot 10^{-14}$ &
25 & 0.1451186573479985220441655$\cdot 10^{-14}$ \\
26 & 0.3888434449455691142432431$\cdot 10^{-15}$ &
27 & 0.1041901454404881356692972$\cdot 10^{-15}$ \\
28 & 0.2791764103483896329393674$\cdot 10^{-16}$
\end{tabular}

As a by-product, insertion of $x=\pm 1$ in (\ref{eq.psiT})
with $\psi(1)=-\gamma$ shows
\begin{eqnarray}
\sum_{k=2}^\infty\frac{1}{k\sqrt{k^2-1}}\frac{k-1+\sqrt{k^2-1}}{k+1+\sqrt{k^2-1}}&=&\frac{1}{2}, \\
\sum_{k=2}^\infty\frac{1}{k\sqrt{k^2-1}}\frac{k+1+\sqrt{k^2-1}}{k-1+\sqrt{k^2-1}}&=&1.
\end{eqnarray}
Linear combinations of these two equations are
\begin{eqnarray}
\sum_{k=2}^\infty
\frac{(k-3)(k+\sqrt{k^2-1})}
{k\sqrt{k^2-1}[(k+\sqrt{k^2-1})^2-1]}
&=&0, \\
\sum_{k=2}^\infty
\frac{k+\sqrt{k^2-1}}
{k\sqrt{k^2-1}[(k+\sqrt{k^2-1})^2-1]}
&=&\frac{1}{8},
\end{eqnarray}
and these two can be combined to
\begin{equation}
\sum_{k=2}^\infty
\frac{k+\sqrt{k^2-1}}
{\sqrt{k^2-1}[(k+\sqrt{k^2-1})^2-1]}
=\frac{3}{8}.
\end{equation}

\bibliographystyle{amsplain}
\bibliography{all.bib}

\end{document}